\documentstyle[12pt]{amsart}
\input epsf

\newtheorem{theorem}{Theorem}[section]
\newtheorem{definition}{Definition}[section]

\newtheorem{lemma}{Lemma}[section]
\newtheorem{proposition}{Proposition}[section]
\newtheorem{corollary}{Corollary}[section]
\newenvironment{proof}{\noindent {\bf Proof:}}{$\Box$}
\newtheorem{conjecture}{Conjecture}[section]

\numberwithin{equation}{section}

\begin{document}

\title{Integral Geometry of
Plane Curves \\ and Knot Invariants }
\author{Xiao-Song Lin}
\address{Department of Mathematics, Columbia University, New York, NY 10027}
\email{xl@@math.columbia.edu}
\author {Zhenghan Wang}
\address {Department of Mathematics, University of Michigan, Ann Arbor,
 MI 48109}
\email{Zhenghan.Wang@@math.lsa.umich.edu}
\date{December, 1994}
\thanks{The first author is supported in part by NSF and the Sloan Foundation
and the second author is supported
by a Rackham Faculty Summer Fellowship at University of Michigan.}

\maketitle

\begin{abstract} We study the integral expression of a knot invariant obtained
as the second coefficient in the perturbative expansion of Witten's
Chern-Simons path integral associated with a knot. One of the integrals
involved turns out to be a generalization of the classical Crofton integral on
convex plane curves and it is related with invariants of generic plane curves
defined by Arnold recently with deep motivations in symplectic and contact
geometry. Quadratic bounds on these plane curve invariants are derived using
their relationship with the knot invariant.
\end{abstract}

\section{Introduction}

The first and second order terms in the perturbative expansion of Witten's
Chern-Simons path integral associated with a knot in the 3-space was first
analyzed by Guadagnini, Martellini and Mintchev [GMM] as well as Bar-Natan [BN]
shortly
after Witten's seminal work. In an announcement appeared in 1992, Kontsevich
perceived a construction of a vast family of knot invariants which,
presumably, contains the same information as the family of coefficients in the
perturbative expansion of the Chern-Simons path integral associated with a
knot [K]. In a recent paper [BT], Bott and Taubes explored this
construction in a much more detailed manner. At this stage,
it seems that a rigorous
foundation has been laid for studying the perturbative expansion of the
Chern-Simons path integral associated with a knot. But, as it seems to us, a
study of each individual knot invariant in this family as concrete and thorough
as possible is what we are lacking of. The first term in the perturbative
expansion turns out to be a classical quantity associated with a space curve
with nowhere vanishing curvature, which was studied extensively under the name
of C\v alug\v areanu-Pohl-White self-linking formula.  Although the
second term as a knot invariant is also classical, we find that the
approach suggested by perturbative theory of the Chern-Simons path integral
provides a deep insight into some of its previously unknown geometric and
topological contents. This will serve as the prototype of our further
investigation.

The knot invariant we study here is, modulo a certain constant, the second
coefficient of the Conway polynomial of a knot. This is a Vassiliev invariant
of second order. So we denote
it by $v_2$. Perturbative expansion of the Chern-Simons path integral leads to
an expression of $v_2$ as the difference of two integrals, $I_X$ and $I_Y$,
over
the knot thought of as a space curve. Or rather, the functional $I_X-I_Y$ on
simple closed space curves can be proved to be a knot invariant and identified
with the second coefficient of the Conway polynomial, modulo a certain
constant.
When the knot approaches to its plane projection, the first integral $I_X$ will
be concentrated at the crossings. On the other hand, the second integral $I_Y$
turns out to be well defined on the plane projection. Actually, it defines a
functional on the space of generic plane curves. Using
the fact that
$I_X-I_Y$ is a knot invariant, we show that $I_Y$ is constant in
each component of the space of generic plane curves. And we can go further to
understand how
$I_Y$ jumps when we pass through the discriminant (non-generic plane curves).
This relates
$I_Y$ and the invariants of plane curves constructed by Arnold [A1,2] with deep
motivations in symplectic and contact geometry.
We believe that this relationship is what Arnold expected in [A2].

We mentioned above that $v_2$ and the second coefficient of Conway polynomial
can be identified only after modulo a certain constant. This constant is the
value of $I_Y$ on a round plane curve. It was first calculated by Guadagnini,
et al. [GMM].  But their computation is lengthy and not illuminating. As we
couldn't understand their computation, we started to look for our own.
It seems that
one should think of $I_Y$ as a
3-dimensional generalization of the Crofton formula (dated 1868) for convex
plane
curves
in integral geometry.
Our computation of $I_Y$ on a round circle is almost parallel to
the classical proof of the Crofton formula.
As the classical proof of the Crofton formula
yields many consequences in integral geometry of
plane curves (e.g., it implies that the
measure of the set of lines intersecting a simple
plane curve is equal to the length of
that curve), we can't stop seeking similar
consequences of our generalized Crofton formula.
Notice that Bott and Taubes [BT] have observed that
the construction of these knot
invariants looks rather similar to the construction
in classical integral geometry. Our
study here seems to make this observation more concrete.

In our investigation of the knot invariant $v_2$, we noticed a quadratic
 bound for the values of $v_2$
on knots with $n$ crossings.
It is
derived from a combinatorial formula for $v_2$.
Such a quadratic bound agrees with
the point of view that Vassiliev invariants should be thought of as
polynomial functions on the set of knots.
We conjecture that a similar bound exists in general. See Section 4.
The combinatorial formula for $v_2$ also leads to quadratic bounds
on Arnold's invariants of plane curves via their
relations with $I_Y$.

The paper is organized as follows. In Section 2, we will define the integral
$I_X$ and $I_Y$ and present some simple calculations and generalizations. In
Section 3, we evaluate the integral $I_Y$ on a round circle in the plane. This
is done by imitating the classical proof of the Crofton formula. In Section 4,
we study the combinatorics of the knot invariant $v_2$ by considering
certain limiting behaviors of the integrals $I_X$
and
$I_Y$. In Section 5, we relate the integral $I_Y$ with
invariants of so called {\em unicursal} plane curves defined
by Arnold [A1].

This work owns its conception to the first author's visit to National Taiwan
University in March, 1994. Discussions with P.-Z. Ong and S.-W. Yang
and their warm hospitality are gratefully acknowledged by the first author.

\section{An integral knot invariant}

For $x\in\Bbb R^3\setminus\{0\}$, we denote by
$$\omega(x)=\frac1{4\pi}\frac{[x,dx,dx]}{|x|^3}$$
the unit area form of the unit 2-sphere $S^2$,
where $[\,\cdot\,,\,\cdot\,,\,\cdot\,]$ is
the mixed product in $\Bbb R^3$.

Let $\gamma:S^1\rightarrow\Bbb R^3$ be a smooth imbedding where we identify
$S^1$ with $\Bbb R/\Bbb Z$. We denote
$$\Delta_4=\left\{(t_1,t_2,t_3,t_4)\,;\,0<t_1<t_2<t_3<t_4<1\right\}$$
and
$$\Delta_3(\gamma)=\left\{(t_1,t_2,t_3,z)\,;\,0<t_1<t_2<t_3<1,\,z\in\Bbb
R^3\setminus\{\gamma(t_1),\gamma(t_2),\gamma(t_3)\}\right\}.$$
Define
\begin{equation}
I_X(\gamma)=\int_{\Delta_4}\omega(\gamma(t_3)-\gamma(t_1))\wedge
\omega(\gamma(t_4)-\gamma(t_2))
\end{equation}

\noindent and
\begin{equation}
I_Y(\gamma)=\int_{\Delta_3(\gamma)}
\omega(z-\gamma(t_1))\wedge\omega(z-\gamma(t_2))\wedge\omega(z-\gamma(t_3)).
\end{equation}

We will call these integrals the $X$-{\it integral} and $Y$-{\it integral}
respectively. They get their names from the
diagrams they correspond to. See Figure 1.

\begin{figure}
\caption{Diagrams for $X$-integral and $Y$-integral.}
\end{figure}

First, we want to simplify the expressions of these two integrals a little
bit .

\begin{lemma} We have
\begin{equation}
\aligned I_X(\gamma)&=-\frac1{(2\pi)^2}\int_{\Delta_4}
\frac{[\gamma(t_3)-\gamma(t_1),\dot\gamma(t_3),
\dot\gamma(t_1)]}{|\gamma(t_3)-\gamma(t_)|^3}\\
&\frac{[\gamma(t_4)-\gamma(t_2),
\dot\gamma(t_4),\dot\gamma(t_2)]}{|\gamma(t_4)-\gamma(t_2)|^3}\,dt_1\,
dt_2\,dt_3\,dt_4.
\endaligned
\end{equation}

\end{lemma}

Notice that $[\gamma(t)-\gamma(t'),\dot\gamma(t),\dot\gamma(t')]$ is
the oriented volume
of the parallelepiped spanned by $\gamma(t)-\gamma(t')$, $\dot\gamma(t)$ and
$\dot\gamma(t')$.

\begin{lemma}
Let
$$E(z,t)=\frac{(z-\gamma(t))\times\dot\gamma(t)}{|z-\gamma(t)|^3}.$$
Then,
\begin{equation}
I_Y(\gamma)=-\frac1{(2\pi)^3}\int_{\Delta_3(\gamma)}\left[E(z,t_1),
E(z,t_2),E(z,t_3)\right]\,d^3z\,dt_1
\,dt_2\,dt_3.
\end{equation}

\end{lemma}

Notice that $E(z,t)dt=dB$ where $B=B(z)$ is the magnetic field induced by
the current
$\gamma(t)$.

Both of these two lemmas come from a straightforward computation.

\setcounter{theorem}{2}
\begin{theorem} Let
$$v_2(\gamma)=I_X(\gamma)-I_Y(\gamma).$$
Then $v_2$ is invariant under isotopy of $\gamma$.
\end{theorem}

So $v_2$ is a knot invariant. This
was proved rigorously by Bar-Natan [BN]
first in his Princeton thesis. See also [BT]. Furthermore, this knot
invariant satisfies
a crossing change formula which identifies itself with the second coefficient
of the Conway polynomial modulo
the constant $v_2(\text{unknot})$. This also justifies the subscription of
$v_2$.

The first step in the proof of this theorem
is to show that both integrals $I_X(\gamma)$ and $I_Y(\gamma)$
are finite. This is done in [BT] by compactifying the
integral domains and showing that the integral forms
extend as smooth forms to the compactification of integral domains. From
this consideration, the following
lemma should be quite obvious.

\addtocounter{lemma}{1}
\begin{lemma} For an immersion $\gamma:S^1\rightarrow\Bbb R^3$ with only
finitely many singularities, the
integral $I_Y(\gamma)$ is finite.
\end{lemma}

\begin{proof} When $\gamma$ is an imbedding, $\Delta_3(\gamma)$ is the total
space of a fibration over $\Delta_3=\{
(t_1,t_2,t_3)\,;\,0<t_1<t_2<t_3<1\}$ whose fibre is $\Bbb R^3$ with three
distinct points deleted. In our case, we
may define $\Delta_3(\gamma)$ similarly and it will have degenerate fibres
over a measure zero set of $\Delta_3$.
Similar to the case where $\gamma$ is an imbedding, $\Delta_3(\gamma)$ with
degenerate fibres deleted can be
compactified and the integral form of $I_Y(\gamma)$ extends to the
compactification. This implies that the integral
$I_Y(\gamma)$ is still finite when $\gamma$ is an immersion with only
finitely many singularities.
\end{proof}

\section{A generalized Crofton formula}

Some remarkable integral formulae associated with convex sets in the plane
are obtained by simple computations of the density $dx\,dy$ in
different coordinate systems.  The classical Crofton formula is
such an example [S].

Let $D$ be a bounded convex set in the plane.  Through each point $P$ exterior
to $D$, there pass two supporting lines of $D$.  Let $s$ and $s'$
respectively be the
lengths of the line segments from $P$ to the corresponding supporting points
$H_{1}$ and $H_{2}$, and
let $\alpha$ be the angle $H_{1}PH_{2}$ between the supporting lines. See
Figure 2. Then
\begin{equation}
\int_{P\notin D} \frac{\sin\alpha}{s\cdot s'}\, dx\,dy=2 {\pi}^{2}
\end{equation}

\begin{figure}
\caption{The setting of the Crofton formula.}
\end{figure}

Let $A=A(P)$ be the area of the parallelogram
spanned by $PH_{2}$ and $PH_{1}$,
then (3.1) can be written as
\begin{equation}
\int_{P\notin D} \frac{A}{s^{2} \cdot {s'}^{2}}\,dx\,dy =2 {\pi}^{2}
\end{equation}

The Crofton formula (3.1) or (3.2) yields many consequences in integral
geometry
of plane curves.  For example, it implied that the measure of the set
of lines intersecting a simple plane curve is equal to the length of
that curve.

It seems amazing that the Crofton integral (left side of (3.1) or (3.2)) is
independent of the shape of $C$. Assuming
this for the moment, let us try to evaluate the Crofton integral when the
boundary of $D$ is the round
circle
$\{(\cos2\pi t,\sin
2\pi t)\,;\,0\leq t\leq1\}$ in the plane.
It is done by a certain change of coordinates.

Let $\phi:\Bbb R^2\setminus
D=\{(x,y)\,;\, x^{2}+y^{2}>1\} \longrightarrow S^{1}\times S^{1} $ be the map
defined by sending each point $P\notin D$ to the
pair of angles $(\theta_1,
\theta_2)$ with $\theta_1<\theta_2$ of the supporting points.
It is easy to check:
\begin{itemize}
\item $\phi$ is one-one;
\item if $P$ goes to infinity in the direction of angle $\theta_0$, then
$\phi (P)$ goes to $(\theta_0 + \pi /2, \theta_0+ 3\pi /2)$;
\item if $P$ approaches the point $(\cos\theta_0, \sin\theta_0)$, then
$\phi (P)$ goes to $(\theta_0, \theta_0)$.
\end{itemize}
Therefore, the image of  $\phi$ covers exactly one half
of the torus $S^1\times S^1$.
By a direct computation, the Crofton integral is the
same as the signed area covered by the image of $\phi$.  It
follows that the Crofton integral is equal to $1/2 \cdot (2\pi)^{2}=
2{\pi}^{2}$.

In general, the same argument will go through since the angle of a
supporting line of the convex set $D$ is well
defined once a center of $D$ is chosen. See [S]

Let $\gamma$ be a simple closed plane curve.  The X-integral $I_{X}$ of
$\gamma$ is $0$ (see (2.3)). The fact that $I_X-I_Y$ is a knot invariant
implies that the $Y$-integral of $\gamma$,
$I_Y(\gamma)$, is invariant under deformation of
$\gamma$ in the space of simple closed plane curves. The computation
of the $Y$-integral for the round circle was first
done by [GMM].
We will provide a computation here
which is similar to the proof of the Crofton formula described above.

Denote $\Bbb R^3_+=\{(z_1,z_2,z_3)\in\Bbb R^3\,;\,z_3>0\}$ and
$\Bbb R^2=\{z_3=0\}\subset\Bbb R^3$.

\begin{figure}
\caption{The setting of the generalized Crofton formula.}
\end{figure}

\begin{theorem} (Generalized Crofton formula)
Let $\gamma =\gamma(t) : S^1\longrightarrow \Bbb R^{2}$
be a simple closed plane curve in $\Bbb R^2$.
Then
\begin{equation}
\int_{\Delta_3}\int_{z\in\Bbb R^3_+}
\frac{V}{\prod_{i=1}^3|z-\gamma(t_i)|^3}\,
d^3z\,dt_1\,dt_2\,dt_3=\frac{{\pi}^3}{6}
\end{equation}

\noindent where $V$ is the oriented volume of the parallelepiped spanned by
$(z-\gamma(t_i))\times\dot\gamma(t_i)$,
$i=1,2,3.$
\end{theorem}

We will call the integral in (3.3) {\it the generalized Crofton integral}.
See Figure 3.

\begin{proof}
Compare with (2.4), the integral in (3.3) is equal to a constant multiple of
$I_Y(\gamma)$. Since
$I_Y(\gamma)=-v_2(\gamma)$ in this case, it is invariant when $\gamma$ is
deformed by an isotopy of the
plane. Therefore, it  suffices to prove the theorem for the round circle
$\{(\cos2\pi t, \sin2\pi t,0)\,;\, 0\leq t \leq 1\} $ in $\Bbb R^2$.

Let $\phi : \Bbb R^{3}_{+}\times(S^1)^3 \longrightarrow
S^{2}_{+}\times S^{2}_{+}\times S^{2}_{+}$ be the map defined by
sending $(z,t_{1}, t_{2}, t_{3})$ to
$$\left(\frac{z-\gamma (t_{1})}{|z-\gamma (t_{1})|},
\frac{z-\gamma (t_{2})}{| z-\gamma (t_{2})| },
\frac{z-\gamma (t_{3})}{| z-\gamma (t_{3})| }\right).$$
Then the generalized Crofton integral on the round circle
is equal to the signed volume of the part of
$(S^{2}_{+})^{3}$ covered by the image of the map
$\phi$ multiplied  by $(2\pi)^3/{3!}$.

{\em Claim:}  There is a subset $A\subset\text{Img}(\phi)$ of full measure in
$(S^2_+)^3$ such that $\phi|\phi^{-1}(A)$ is one-one.

It follows from the claim that the generalized Crofton integral on the round
circle
is $\pm (2\pi)^3/3!\cdot (1/2)^{3}=\pm {\pi}^3/6$.
By checking the orientations, we know the integral equals ${\pi}^3/6$.

{\em Proof of the claim:}  Let $v$ be a vector in the upper hemi-sphere
$S^{2}_{+}\subset\Bbb R^{3}$, and $\phi_1:\Bbb R^{3}_{+}\times S^1
\longrightarrow S^{2}_{+}$ be the map defined by sending
$(z,t)$ to
$$\frac{z-\gamma (t)}{| z-\gamma (t)|}.$$
Then ${\phi_1}^{-1}(v)$ is the half infinite cylinder
$$C_{v}=\{\gamma (t) + sv\,;\, 0\leq t\leq 1,\,s\geq0\}.$$
Let $(v_{1}, v_{2}, v_{3})$ be a point in $(S^{2}_{+})^{3}$.  Then
$\phi^{-1}(v_{1}, v_{2}, v_{3})$ is in one-one correspondence with the set
of intersections
of the three half infinite cylinders $C_{v_{1}}$, $C_{v_{2}}$, and $C_{v_{3}}$.
If $v_{1} \neq v_{2}$, then $C_{v_{1}}$, $C_{v_{2}}$ intersect
in an arc lying on both $C_{v_1}$ and $C_{v_2}$ whose ends are a pair of
antipode points of the round
circle $\gamma$.  If $v_{1}$, $v_{2}$ and $v_{3}$ are pairwise
distinct, then $C_{v_{1}}$ and $C_{v_{2}}$ intersect in an arc $A_{12}$
on $C_{v_{1}}$, and similarly, $C_{v_{1}}$ and $C_{v_{3}}$ intersect
in an arc $A_{13}$ on $C_{v_{1}}$. Let $p:\Bbb R^3\rightarrow\Bbb R^2$ be the
projection.
If $p(v_2)-p(v_1)$ and $p(v_3)-p(v_1)$ are not
colinear ($p(v_1)$, $p(v_2)$ and $p(v_3)$ are in general position), then
the ends of
$A_{12}$ and
$A_{13}$ form two distinct pairs of antipode points on the round circle.
This implies that
$A_{12}$ and
$A_{13}$ intersect at exactly one point in $R^{3}_{+}$.
So if $p(v_1)$, $p(v_2)$ and $p(v_3)$ are in general position, then
$\phi^{-1} (v_{1}, v_{2}, v_{3})$ is a single point in
$\Bbb R^3_+\times(S^1)^3$.
The remnant in $(S^{2}_{+})^{3}$ is of measure 0.
This proves the claim and thus the theorem.
\end{proof}

It seems very likely that the argument in the
proof of the claim above also works if $\gamma$ is a convex curve in $\Bbb
R^2$. If this is so, we will have a direct proof of
the generalized Crofton formula for convex plane curves. Here is a
very interesting intuitive interpretation of the claim
in the proof of Theorem 3.1. Image that a fixed round circle in
the plane starts to move in the plane with
three non-colinear constant velocities respectively so that we will see
three round circles of the same radius in a moment.
Then there will be exactly one moment when these three circles have
exactly one intersection. Intuitively, this
should also be true if we start with a convex curve. We are not
sure whether such a phenomenon has been discussed in the literatures.

\setcounter{corollary}{1}
\begin{corollary} $v_2(\text{\em unknot})=-1/24$.
\end{corollary}

\begin{proof} It is easy to see that on a round circle,
$I_Y$ is equal to the generalized Crofton integral on the
round circle times
$2\cdot1/(2\pi)^3=1/4\pi^3$. So $I_Y$ on a round circle is 1/24.
\end{proof}

\section{The combinatorics of the integral knot invariant $v_{2}$}

As we mentioned before, the knot invariant $v_2$ can be
identified, modulo the constant $v_2(\text{unknot})=-1/24$,
with the second coefficient of the Conway polynomial via
a crossing change formula. From this identification, one can
draw most of the conclusions about $v_2$ in this section.
Therefore, the main interest of this section
is probably to see how one can study $v_2$ by studying
certain limiting behaviors of the integrals $I_X$ and $I_Y$.
Such a consideration also leads to the discovery of
the relationship between the knot invariant $v_2$ and invariants
of plane curves discussed in the next section.

Let $\gamma : S^{1} \longrightarrow\Bbb R^{3}$ be an imbedding.
Since it is very difficult
to compute $v_{2}$ by evaluating both integrals
$I_X(\gamma)$ and $I_Y(\gamma)$ directly,
we study a limiting situation when the
curve is pushed into a plane via a regular projection.
It turns out that the limits of both integrals can be computed, and
this gives us a new way of studying the knot invariant $v_{2}$.
The same analysis applied to the Gauss linking formula
give us the well-known combinatorial formula
for the linking number.

When an imbedded curve $\gamma$ in $\Bbb R^{3}$ acquires a double point, the
$X$-integral of $\gamma$ blows up.  On the other hand, the $Y$-integral
is still
meaningful by Lemma 2.4.

\begin{proposition}
Let $\gamma : S^{1}\longrightarrow\Bbb R^{3}$ be an immersion
with only transverse
double points.  Then

(1)  $I_{Y}(\gamma)$ does not depend on the parameterization or
the orientation of $\gamma$.

(2) Let $T: \Bbb R^{3}\longrightarrow\Bbb R^{3}$ be an affine similarity.
Then $I_{Y}(T\circ \gamma)=I_{Y}(\gamma)$.
\end{proposition}

The proof is immediate.

Let $K$ be the knot type of the imbedding $\gamma:S^1\longrightarrow
\Bbb R^3$. Thought of as a knot diagram,
$K$ can be drawn inside the plane except around each crossing.
Assume that when we make an over-pass at a crossing
of $K$, we go along a semi-circle of radius $\epsilon$ perpendicular to the
plane. The other parts of $K$ lie
completely in the plane. Denote such a diagram of $K$ by $K_{\epsilon}$.
Also, we associate to each crossing a sign $\pm 1$
by the usual right-hand-rule.

When $\epsilon$
approaches 0, $K_{\epsilon}$ limits to a closed plane curve with only
transverse
double points.   Denote this limiting plane curve by $K_{0}$.
We associate to each double point
the sign of the corresponding crossing, and we call $K_{0}$
with these signs the {\em signed limiting plane curve.}
In general, let $C$ be an immersed circle in the plane with only
transverse double points.  If each double point is associated
with a sign, then $C$ is called {\em a signed immersed circle}.
For each closed plane curve with only transverse double points, we have
a chord diagram defined as follows.

\setcounter{definition}{1}
\begin{definition}
Let $C$ be an immersed circle with only transverse double points
in the plane.  The {\em chord diagram} (or {\em Gauss diagram}) of $C$
is the combinatorial pattern of a finite
collection of chords with both ends sticking to the
circle, connecting the points
sent by the immersion to the same double point of the immersed
circle.
If each chord is associated with a sign, the chord
diagram  is called a
 {\em signed chord diagram.}
\end{definition}

If $K_{0}$ is the limiting plane curve of a knot diagram $K_\epsilon$, then we
associate to each chord of the chord
diagram of $K_{0}$ the sign of the corresponding double
point.
For each pair of chords across to each other in the signed
chord diagram of $K_{0}$,
we assign to this pair of chords a sign equal to the product of the signs
of the two chords.

\addtocounter{proposition}{1}
\begin{proposition}
Let $K$ be a knot diagram with n crossings.  Then
the limit $$\lim_{\epsilon \rightarrow 0} I_{X}(K_{\epsilon})$$
exists.  Let $I_{X}(K_{0})$ denote this limit.  Then
\begin{equation}
I_{X}(K_{0})=\frac{n}{16} + \frac{(c_{+} -c_{-})}{4}
\end{equation}

\noindent where $c_{+}$ ($c_-$, respectively) is the
number of pairs of chords in the
signed chord diagram of $K_{0}$ with a positive (negative, respectively) sign.
\end{proposition}

\begin{proof} Our first observation is that when $\dot\gamma(t)$ and
$\dot\gamma(t')$
are co-planar, then the Gaussian form $\omega(\gamma(t)-\gamma(t'))$
is 0.

We may assume that the $i$-th crossing of $K_{\epsilon}$ looks like the
crossing depicted in Figure 4. Let
$C_{\epsilon}^{i}$ be the semi-circle of radius $\epsilon$ at that crossing
and $A_a^{i}$ be
a line segment
$[-a, a]$ with a fixed small $a>0$ on $K_{\epsilon}$ running under
$C_{\epsilon}^{i}$. Furthermore, denote by
$C_{\epsilon,a}^{i}$ the union of $C_{\epsilon}^{i}$ and line segments
$[-a,-a+\epsilon]$ and $[a-\epsilon,a]$.

\begin{figure}
\caption{The local picture of a knot projection at a crossing.}
\end{figure}

By the observation above, if both $\gamma(t)$, $\gamma(t')$ lie outside
$\cup_{i=1}^{n} C_{\epsilon}^{i}$, or
they both lie inside some $C_{\epsilon,a}^{i}$,
then $\omega (\gamma(t)-\gamma(t'))=0$.  If $\gamma(t) \in C_{\epsilon}^{i}$,
and $\gamma(t')$ lies outside
$C_{\epsilon,a}^{i} \cup A_a^{i}$, then $|\gamma(t) -
\gamma(t')| $ is bounded from below by a constant.  It follows that when
$\epsilon$ approaches $0$, the integral $I_{X}$ over all those pairs
goes to $0$.
Note that we need to fix a base point in order to evaluate
the $X$-integral $I_{X}$ for a curve.  So we choose a base point on
$K_{\epsilon}$ which is not in any $C_{\epsilon,a}^{i}$ or
$A_{a}^{i}$.
Then all nonzero limits come from the following two
cases:

(1)  $\gamma(t_{1}), \gamma(t_{2}) \in C_{\epsilon}^{i}$, and
$ \gamma(t_{3}), \gamma(t_{4}) \in A_a^{i}$;

(2) $\gamma(t_{1})\in C_{\epsilon}^{i}, \gamma(t_{3})\in A_a^{i}$ or vice
versa,
and $\gamma(t_{2})\in C_{\epsilon}^{j}, \gamma(t_{4})\in A_a^{j}$ or vice
versa, with $i\neq j$.

By a direct computation, the limit for case (1) is always $1/16$.
It is independent of the sign of the crossing.

Since $t_{1}< t_{2} <t_{3} < t_{4}$,  case (2) is possible if and only
if the chords corresponding to the $i$ and $j$ cross each other.
In this case, the limit is $\epsilon_{i} \epsilon_{j}/4$,
where $\epsilon_{i}, \epsilon_{j}$ are the signs of the $i$-th and $j$-th
crossings.  This completes the proof.
\end{proof}

Let $C$ be an immersed circle in the plane
with only transverse double points.  Then we can resolve $C$ to knots by
changing each double point to a crossing.  There are $2^{n}$
resolutions if $C$ has $n$ double points.

\setcounter{corollary}{3}
\begin{corollary}
Let $C$ be an immersed circle in the plane with only transverse
double points.

(1) If $C$ is resolved to two knots $K^{1}$ and $K^{2}$ of the
same knot type,  then $I_{X}(K^{1}_{\epsilon})$ and
$I_{X}(K^{2}_{\epsilon})$
have the same limit.

(2) $I_Y(C)$ is invariant when $C$ is deformed in the plane without changing
its chord diagram.

(3) If $C$ is resolved to
an unknot $K^u$, then
$$I_{Y}(C)=\frac{1}{24}+ \frac{n}{16} + \frac{(c_{+}- c_{-})}{4}$$
where $c_{\pm}$ are computed using the signed chord diagram of $K^u$.
\end{corollary}

Corollary 4.4 (2) implies thousands
of integral formulae like the generalized Crofton formula.

\begin{proof} (1) This is because the limit of $I_Y(K_\epsilon^i)$ when
$\epsilon$ approaches to 0
is $I_Y(C)$, for $i=1,2$, and $v_2(K^1)=v_2(K^2)$.

(2) Resolve $C$ to a knot $K$, since both $I_X(K_0)$ and $v_2(K)$ are
 invariant when $C$ is deformed
in the plane
without changing its chord diagram, so does
$I_Y(C)=I_X(K_0)-v_2(K)$.

(3) This is a direct consequence of Corollary 3.2 and  Lemma 4.3.
\end{proof}

Let $K$ be a knot diagram. At each crossing of $K$, the modification of
the knot diagram depicted in Figure 5
changes $K$ into a link of two oriented components. Such a modification
is called an {\it oriented surgery} at a
crossing.

\begin{figure}
\caption{An oriented surgery.}
\end{figure}

\begin{corollary}
Let $K$ be a knot diagram with n crossings. And $K^{u}$ be an unknot
obtained from $K$ by changing some crossings.

(1) Let $l_{i}$ be the linking number of the two component link
obtained from the oriented surgery at the i-th crossing of $K$, and
$l_{i}^{u}$ be the corresponding linking number from $K^{u}$.  Then
\begin{equation}
v_{2}(K)=\frac{1}{2} \sum_{i=1}^{n} (l_{i} -l_{i}^{u})-\frac{1}{24};
\end{equation}

(2) $|v_{2}(K)| \leq n(n-1)/4 +1/24;$

(3) $v_{2}(K)+1/24$ is an integer;

(4) $v_2(K)$ is independent of the orientation of $K$.

\end{corollary}

\begin{proof}
Let $K_{0}$ be the limiting plane curve of $K$.
And $K_{0}^{u}$ be the limit plane curve of $K^{u}$.
By Lemma 4.3 and Corollary 4.4 (3),
\begin{equation}v_{2}(K)=\frac{c_{+}(K_{0}) -c_{-}(K_{0})}{4} -
\frac{c_{+}(K_{0}^{u}) -c_{-}(K_{0}^{u})}{4}-\frac{1}{24}.
\end{equation}

The difference between the signed chord diagrams of
$K_{0}$ and $K_{0}^{u}$ is that some signs of chords are changed.

(2)  As the chord diagrams of $K_{0}$ and $K_{0}^{u}$ both
have n-chords, there are at most
${n(n-1)}/2$ intersections among chords.
Thus, (4.3) implies
$| v_{2}(K) | \leq {n(n-1)}/{4}+1/24.$

(3) Let $d_{i}$ be the number of intersections of all chords
with the $i$-th chord counted with signs.
Then $d_{i}=2l_{i}$, so $d_{i}$ is always an even integer.

{\em Claim:}  Let $K_{0}^{1}$ be a signed plane curve obtained
from $K_{0}$ by changing one sign of the double points.
Then $I_{X}(K_{0})-I_{X}(K_{0}^{1})$ is an integer.

{\em Proof of the claim:}  By (4.1), we need only count
the changes of the intersections between the chord corresponding to the
double point where the sign is changed
and other chords, divided by 4.  When the sign of the $i$-th
chord is changed, $d_{i}$ is changed to $-d_{i}$.  As
$d_{i}$ is even,  the total change $2d_{i}$ is divisible by 4.

Now we finish the proof of (3) as follows:
choose a sequence of double points so that when we change
the their signs one after another to get a sequence of
signed plane curves $K_{0}^{1}, \cdots , K_{0}^{s}$ with
$K_0=K_0^1$ and $K_0^s$ is the limit of an unknot. Then
$$v_{2}(K)+\frac1{24}=I_{X}(K_{0}) -I_{X}(K_{0}^{s})=
\sum_{i=1}^{s-1}\left(I_X(K^i_0)-I_X(K_0^{i+1})\right).$$
Thus the
proof of (3) is  completed.

Finally, it is clear that (1) follows from the proof of (3).
It is also clear that (4) can be derived in many ways and
one way is via the formula (4.2) since the linking
number will not change if one changes the orientations of both
components of a link.
\end{proof}

Corollary 4.5 (2) is of particular interest to us. Recall
that a knot invariant is of finite type or a Vassiliev
invariant if it vanishes on \lq\lq higher order differences of knots".
One may think of such a knot invariant as a
\lq\lq polynomial function" on the set of knot types. See, for example,
[B]. Since $v_2$ is known to be a
Vassiliev invariant of order 2, the bound for its values
on knots with $n$ crossings agrees with such a point of
view. Notice that using the combinatorial formula for a Vassiliev invariant
of order 3 given by Lannes [L], we
can get a similar bound for values of Vassiliev invariants of
order 3 on knots with $n$ crossings. This leads us to
the following conjecture.

\setcounter{conjecture}{5}
\begin{conjecture} For every Vassiliev invariant
of order $k$, say $v_{k}$, there is a constant $C$ such that if $K$ is a knot
with $n$ crossings, then \[
| v_{k}(K) | < C n^{k}.\]
\end{conjecture}

Combinatorial formulae for Vassiliev invariants of
lower orders were also discussed by Polyak and Viro.

\section{Invariants of unicursal curves}

Considering the space $\cal M$ of all immersions of
$S^1$ into the plane. By a classical theorem of Whitney,
components of this space can be indexed by $\Bbb Z$ using
the winding number. We will denote by $\cal M_w$ the
component of $\cal M$ whose members all have winding number or index $w$.
If we want to look at
$\cal M$ more carefully, we will see generic immersions and non-generic
ones. Generic immersions are those with only
transverse double points. The set $\Sigma$ of
non-generic immersions or the discriminant of $\cal M$
can be thought of as a stratified space whose top stratum have
three components of particular interest to us. One component
consists of immersions with exactly one transverse triple
point, and the other two consist of, respectively, immersions
with exactly one direct or inverse self-tangency point
where two  tangent branches of the curve concave in opposite direction.
A self-tangency point is called {\it direct},
if two tangent vectors of the curve at the tangent point are in the same
direction, or it is called
{\it inverse} otherwise. It turns out that these three
components of the top stratum of
$\Sigma$ corresponding respectively to immersions
with exactly one transverse triple point, one direct self-tangency
point or one inverse self-tangency point are all well \lq\lq co-oriented".
This means that we can talk about the
positive or negative side of these three component of the top stratum of
$\Sigma$ in $\cal M$. It is quite easy to see that a path in $\cal M$ can be
perturbed so that it only crosses
$\Sigma$ transversally through these three
components at finitely many places. For a detailed study of the topology of
$\Sigma\subset\cal M$, see [A1].

We will call a generic immersion of $S^1$ into the plane a {\it unicursal
curve}. Two unicursal curves are equivalent
if they belong to the same component of $\cal M\setminus\Sigma$. It is not
hard to see that a path in $\cal
M\setminus\Sigma$ is the same as a deformation of a unicursal curve without
changing its chord diagram. An
{\it invariant of unicursal curves} assigns values to every unicursal curve
and equivalent unicursal curves
should be assigned with the same value.

\begin{lemma} $C\rightarrow I_Y(C)$ is an invariant of unicursal curves.
\end{lemma}

This is simply a restatement of Corollary 4.4 (2).

In [A1], Arnold constructed three basic invariants of unicursal curve, $St$ and
 $J^{\pm}$. Up to an additive
constant, they are completely determined
by the way they jump when a deformation of unicursal curves crosses through a
triple point, or a direct or inverse
self-tangency point. To describe these invariants, we need to make some
definition first.

\setcounter{definition}{1}
\begin{definition} (1) A transversal crossing of a self-tangency point is {\em
positive} if the number of double points
grows (by 2).

(2) A transversal crossing of a triple point is {\em positive} if the new-born
vanishing triangle is positive.
\end{definition}

Here for a given a unicursal curve $C$, a {\it vanishing} triangle of $C$ is a
triangle formed by three branches of $C$ and no other branches of $C$ are
allowed to run into such a triangle. At a
transversal crossing of a triple point, one sees the death of one vanishing
triangle and the birth of another one.  The
sign of a vanishing triangle is defined as follows.
The orientation of the immersed circle defines
a cyclic ordering of the sides of the vanishing triangle.
Hence the sides of the triangle acquire orientations induced by the
ordering.  But each side has also its own direction which might
coincide, or not, with the orientation defined by the ordering.
For each vanishing triangle, let $q$ be the number of sides
equally oriented by the ordering and their directions.  Then
the {\em sign} of a vanishing triangle is $(-1)^{q}$. It is easy to check that
at a transversal crossing of a triple
point, the dying vanishing triangle and the new-born vanishing triangle always
have opposite signs.

\setcounter{theorem}{2}
\begin{theorem} (Arnold) (1) There exists a unique (up
to an additive constant)
invariant of unicursal curves of
fixed index whose value remains unchanged at a transversal crossing of a
self-tangency point, but increases by 1 at a
positive transversal crossing of a triple point. This invariant is denoted by
$St$ with an appropriate normalization.

(2) There exists a unique (up to an additive constant) invariant of unicursal
curves of fixed index whose value
remains unchanged at a transversal crossing of a triple point or an inverse
(respectively, direct) self-tangency point,
but increases by 2 (respectively, $-2$) at a positive transversal crossing of
a direct (respectively, inverse)
self-tangency point. This invariant is denoted by $J^+$ (respectively, $J^-$)
with an appropriate normalization. The
invariants $J^+$ and $J^-$ are related by $J^+-J^-=n$ on unicursal curves with
$n$ double points.

(3) These invariants are independent of the orientation of unicursal curves.
\end{theorem}

Here a normalization means to choose a unicursal curve $C_w$ for each index
$w\in\Bbb Z$ and the value of
the invariant in question on $C_w$. One may think of these three invariants
$St$ and $J^{\pm}$ of unicursal curves as
dual to those three components of the top stratum of $\Sigma$ corresponding to
one triple point, one direct
self-tangency point and one inverse self-tangency point respectively. Any
invariant of unicursal curves which jumps by a
constant at a transversal crossing of a triple point and a self-tangency point
can be expressed uniquely as a linear
combination of $St$ and $J^{\pm}$, modulo a constant depending on the index.

To simplify the terminology, we define several operations on unicursal curves
analogous to the Reidermeister
moves in knot theory. See Figure 6. A {\em type} I {\em move} on a unicursal
curve kills one small kink on it.
A {\em type} II$^+$ (II$^-$, respectively) {\em move} is a positive
transversal
crossing of a direct (inverse,
respectively) self-tangency point. Finally, A {\em type} III {\em move} is a
positive transversal crossing of a
triple point.

\addtocounter{definition}{1}
\begin{definition}
Let $C$ be a unicursal curve.  Then we define $$\alpha(C)=I_{Y}(C) + \frac{n}
{16}- \frac{1}{24}$$
where $n$ is the number of double points of $C$.
\end{definition}

As both the $Y$-integral and the number of double points are
invariants of unicursal curves, so is $\alpha$.
If $C$ is resolved to an unknot $K^u$, then
\begin{equation} \alpha(C)=
\frac{n}{8} + \frac{c_{+} -c_{-}}{4}
\end{equation}

\noindent where $c_{\pm}$ are computed using the signed chord diagram
of $K_0^u$.
{}From this formula, we see that each double point of $C$ contributes 1/8 to
$\alpha$,
and each pair of intersecting chords contributes $\pm 1/4$ to $\alpha$.

\begin{figure}
\caption{Elementary moves among unicursal curves.}
\end{figure}

\addtocounter{theorem}{1}
\begin{theorem}
The invariant $\alpha$ of unicursal curves has the following properties:

(1) $\alpha$ equals to 0 for every simple closed plane curve;

(2)  $\alpha$ is decreased by $1/8$ if a type {\em I }
move is performed;

(3)  $\alpha$ is unchanged if a type {\em II$^+$ }
move is performed;

(4)  $\alpha$ is decreased by 1/4 if a type {\em III}
move is performed.

As a consequence of (2), (3), and (4), $\alpha$
is increased by $1/4$ if a type {\em II$^-$ }
move is performed. Furthermore, we have

(5) $|\alpha(C)|\leq n^2/8$ where $n$ is the number of double points on $C$;

(6) $\alpha(C)$ is independent of the orientation of $C$.
\end{theorem}

\begin{proof}
They are all consequences of (5.1).

(1) This is obvious.

(2) Note that the chord corresponding to the double point
in a type I move is an isolated chord (it does not intersect
any other chord).  So this double point contributes only $1/8$
to $\alpha$.

(3) We may assume that after a type II$^+$ move, the two new double points are
resolved with opposite signs. They
contribute 1/4 to $\alpha$. But the two chords of the new double points
intersect and every other chord either
intersects them both or misses them both. So their contribution to $\alpha$
is -1/4. This implies (3)

(4) Assume that a type III move changes $C$ to $C'$ and the vanishing
triangles on $C$ and $C'$ are resolved as
depicted in Figure 7. This is done by choosing a base point and resolving
double points according to the rule that the branch one walks through first
is always
above the branch one walks through second. The resulting knot is an unknot.

On the level of chord diagrams, there are two cases to study. See again
Figure 7.
In both cases, the chord diagrams for $C$ and $C'$ have the same number of
chords.
Let
$\{a,b,c\}$ and
$\{a',b',c'\}$ be the signed chords at the vertices of the vanishing triangles
of $C$ and $C'$ respectively. The edges
of the vanishing triangle of $C$ (respectively, $C'$) correspond to three
disjoint arcs on $S^1$ and each chord in
$\{a,b,c\}$ (respectively, $\{a',b',c'\}$) connects end points of two distinct
arcs. No other chords of $C$
(respectively, $C'$) with touch these arcs. Furthermore, when a type III move
changes $C$ to $C'$, the chord diagram
of $C$ is changed to the chord diagram of $C'$ by switching every two end
points of $\{a,b,c\}$ paired as the end
points of those arcs on $S^1$. The signs of chords will not be changed. So
the contribution of $\{a,b,c\}$ to
$\alpha(C)$ will be 1/4 more than the contribution of $\{a',b',c'\}$ to
$\alpha(C')$.  Furthermore, if another chord
intersects one chord in $\{a,b,c\}$, it will also intersects one in
$\{a',b',c'\}$ with the same sign. Therefor, $\alpha(C')=\alpha(C)-1/4$.

(5) This is a consequence of (5.1) and the bound $|c_+-c_-|\leq n(n-1)/2$.

(6) This also follows easily from (5.1).  \end{proof}

\begin{figure}
\caption{The proof of Theorem 5.5 (4).}
\end{figure}

Now by Theorem 5.4, it is very easy to compute $\alpha$.  Consequently, it is
very easy to compute $I_{Y}(C)$ for any unicursal curve $C$.
For example, $\alpha(\infty)=1/8$.  It follows that
$I_{Y}(\infty)=5/48$. Here the symbol $\infty$ is used to denote a unicursal
curve of the same shape.

\setcounter{corollary}{5}
\begin{corollary} We have
\begin{equation} \alpha =- \frac{2St + J^{-}}{8}.
\end{equation}

\end{corollary}

\begin{proof} The invariant in the right side of (5.2) changes in the same
way as $\alpha$ does under type II$^{\pm}$ and III moves. Checking the
initial values of $St$ and $J^{\pm}$ given in [A1]
verifies (5.2).
\end{proof}

A result of F. Aicardi (see [A1]) says that $2St+J^+=0$ holds if the chord
diagram
of $C$ have no intersecting chords. Here is a generalization of this result.

\begin{corollary} The identity $2St+J^+=0$ holds for a unicursal
curve $C$ with $n$ double points if and only if
$\alpha(C)=n/8$, and if and only if a certain signed chord diagram coming
from an unknot resolution of $C$ have
$c_+=c_-$.
\end{corollary}

\begin{proof} For unicursal curves with $n$ double points, we have
$J^+-J^-=n$. So (5.2) gives us
\begin{equation}
\alpha =- \frac{2St + J^{+}}{8}+\frac{n}{8}.
\end{equation}

\noindent on unicursal curves with $n$ double points. This together with
(5.1) proves the corollary.
\end{proof}

Arnold's triple $St,\,J^+,\,J^-$ are related via some other
index-type invariants. They are defined as follows.

Let $C$ be a unicursal curve such that at each transverse double point,
the two tangent vectors of $C$ are orthogonal.
At every double point, we may divide $C$ into two branches
$C_1$ and
$C_2$. In fact, the preimage of that double point on $S^1$ cuts $S^1$ into
two arcs and the images of these two arcs
under the immersion are $C_1$ and $C_2$
respectively. These two branches $C_1$ and $C_2$ are ordered such that if
the  outgoing tangent vectors of
$C_1$ and
$C_2$ at the corresponding double point are $v_1$ and $v_2$, respectively,
then the frame $\{v_1,v_2\}$ has the
same orientation as that of the plane.

\addtocounter{definition}{3}
\begin{definition}
The {\em half-index} $i_{1}$ (respectively, $i_{2}$) of a double point
is the angle of the rotation of the radius-vector connecting the double point
to a point moving along $C_1$ (respectively, $C_2$) from the double
point to itself divided by $\pi/2$.  The {\em index} of a double
point is the difference $i=i_{1}-i_{2}$.
\end{definition}

The invariants $I^{\pm}$ are defined to be
\begin{equation}I^{\pm}=\frac{\sum i \pm 2n}{4}\end{equation}
where $n$ is the number of double points, and the sum is over
all double points.

Note that $I^{+}-I^{-}=n$. And we have
\begin{equation}
J^{\pm}=I^{\pm}-3St
\end{equation}
as shown in [A1]. So, among these three invariant $St$, $J^+$ and $J^-$,
there is only one which is essentially not of
index-type. Corollary 5.10 below shows that  $\alpha$ is such an invariant.

The following theorem is essentially from [A1].

\addtocounter{theorem}{3}
\begin{theorem}
$I^{-}$ is determined by the following properties:

(1) $I^{-}$ of a simple closed plane curve is 0;

(2) $I^{-}$ is decreased by $(i-2)/4$
if a type {\em I} move is performed, where $i$ is the index of this double
point;

(3) $I^{-}$ is unchanged if a type {\em II$^+$} move is performed,
and $I^{-}$ is decreased by 2 if a type {\em II$^-$} move is performed;

(4) $I^{-}$ is increased by 3 if a type {\em III} move is performed.

\end{theorem}

\addtocounter{corollary}{2}
\begin{corollary}
We have the following identities:

(1) $St= I^{-} +8\alpha $;

(2) $J^{-} = -2I^{-} -24\alpha$;

(3) $J^{+} =n -2I^{-} -24\alpha$.

\end{corollary}

\begin{proof}
The proof follows from (5.2), (5.5) and $J^+-J^-=n$.
\end{proof}

Let $C$ be a unicursal curve, and $x$ be a point which is not
on $C$.  The {\em winding number of $C$ relative to $x$}
is the degree of the position map
$$S^1\longrightarrow S^1:\,\, t\mapsto\frac{C(t)-x}{|C(t)-x|}.$$
The relative winding number remains unchanged if $x$ moves in a connected
component of $\Bbb R^{2}\setminus C$.
We will use the relative winding number to estimate the index of double
points.

Let us notice that the oriented surgery at a crossing on a knot diagram can
be generalized to the oriented surgery at
a double point on a unicursal curve. The oriented surgery at a double point
$x$ on a unicursal curve $C$ will result
two new unicursal curves $C_1$ and $C_2$ intersecting each other
transversally. The component of $\Bbb
R^2\setminus C_1\cup C_2$ where $x$ lies is well defined.

\addtocounter{lemma}{9}
\begin{lemma}
(1) Let $w_{1}$ and $w_{2}$
be the winding numbers of $C_{1}$ and $C_{2}$ relative to $x$.  Then
the index of this double point $x$ on $C$ is
$4w_{1} -4w_{2} -2$.

(2) The index $i$ of any double point
on a unicursal curve with $n$ double points satisfies the inequality
$|i|\leq 4n+6$.
\end{lemma}

\begin{proof}
(1) If $$i_1=\frac{\theta_1}{
\pi/2}\qquad\text{and}\qquad i_2=\frac{\theta_2}{\pi/2},$$ then
$$w_1=\frac{\theta_1-\pi/2}{2\pi}
\qquad\text{and}\qquad w_2=\frac{\theta_2+\pi/2}{2\pi}.$$
So we have $i=i_1-i_2=4w_1-4w_2-2.$

(2) Let $C$ be a unicursal curve.
Then $C$  divides the plane into many regions.
Inside each region, pick a point.  Then we assign
to each region the winding number
of $C$ relative to this point.

{\em Claim:}  If $C$ has $n$ double points, then
the maximum of the absolute value
of the winding numbers for all regions is no greater than $n+1$.

{\em Proof of the claim:}  It is easy to check that
if two regions are adjacent with a common edge,
then their winding numbers differ by 1.  If $C$
has $n$ double points, then $C$ has $2n$ edges, and
the plane is divided into $n+2$ regions (including the unbounded one).
Note that the winding number for the unbounded region
is always 0, and the claim follows.

Now the inequality in (2) follows from (1) and the claim.
\end{proof}

\addtocounter{corollary}{1}
\begin{corollary} We have the following bounds on unicursal curves with $n$
double points:

(1) $| I^{-} | \leq n^{2}+2n $;

(2) $ | St | \leq 2n^2+2n$;

(3) $| J^{-} | \leq 5n^2+4n$.

\end{corollary}

These follow easily from Theorem 5.5 (5), Corollary 5.10 and Lemma 5.11.

It is an interesting question to study the extremal curves
for each invariant.
We can fix either the number of double points, or fix the
{\em index} of the curve.  For some conjectures about
upper bounds of $St, J^{\pm}$,
see [A1].  The extremal values of $\alpha$ are still unknown.  The
examples in Figure 8 show that there are curves with $\alpha \leq 0$.
Actually, there is no bound from below for $\alpha$ if we do not
fix the number of double points.  Notice that both
curves have only positive triangles.  There are curves with only negative
triangles, too.

\begin{figure}
\caption{Examples of unicursal curves with $\alpha\leq0$.}
\end{figure}

\end{document}